\documentclass[12pt]{article}
\usepackage{amsmath}
\usepackage{amssymb}
\usepackage{cite}
\newsymbol \blackbox 1004
\newcommand{\eh}{\hfill}\newlength{\sperr}

\newenvironment{proof}{{\settowidth{\sperr}{\bf\rm
Proof}%
\par\addvspace{0.3cm}\noindent\parbox[t]{1.3\sperr}
{\bf\rm P\eh r\eh o\eh o\eh f\eh }%
}}{\nopagebreak\mbox{}
$\blackbox$\par\addvspace{0.3cm}}

\def\nn{\nonumber}
\def\a{\alpha}
\def\b{\beta}

\def\Lam{\Lambda}
\def\s{\sigma}
\def\la{\lambda}
\def\om{\omega}

\def\S{\Sigma}
\def\t{\theta}
\def\Up{\Upsilon}
\def\vp{\varphi}
\def\vt{\vartheta}
\def\ve{\varepsilon}
\def\wh{\widehat}

\def\ov{\overline}

\def\BC{{\mathbb C}}
\def\BR{{\mathbb R}}
\def\BN{{\mathbb N}}
\def\clp{{\mathcal P}}
\def\cla{{\mathcal A}}

\def\clh{{\mathcal H}}

\def\cln{{\mathcal N}}

\def\cld{{\mathcal D}}

\def\im{{\rm Im\ }}

\newcommand{\E}{\mathrm{e}}
\newcommand{\I}{\mathrm{i}}
\def\mf{\mathfrak}

\newtheorem{Pa}{Paper}[section]
\newtheorem{Tm}[Pa]{{\bf Theorem}}

\newtheorem{Cy}[Pa]{{\bf Corollary}}
\newtheorem{Rk}[Pa]{{\bf Remark}}

\newtheorem{Dn}[Pa]{{\bf Definition}}
\newtheorem{Pn}[Pa]{{\bf Proposition}}

\title{Weyl theory and explicit solutions of direct and inverse
problems for a Dirac system with rectangular matrix potential  }

\author{B. Fritzsche, B. Kirstein, I.Ya. Roitberg, A.L. Sakhnovich}

\date{}
\begin{document}
\maketitle

\begin{abstract} A non-classical Weyl theory is developed for
Dirac systems with rectangular matrix potentials. The notion
of the Weyl function is introduced and the corresponding
direct problem is treated. Furthermore, explicit solutions of the 
direct and inverse problems are obtained for the case of rational
Weyl matrix functions.

\end{abstract}

{MSC(2010):}  4B20, 34L40, 15A15, 93B15.


\section{Introduction} \label{intro}
\setcounter{equation}{0}
Consider  
self-adjoint Dirac-type (also called Dirac, ZS or AKNS)
system, which is
a classical matrix differential equation:
\begin{equation}       \label{1.1}
\frac{d}{dx}y(x, z )=\I (z j+jV(x))y(x,
z ) \quad
(x \geq 0),
\end{equation}
\begin{equation}   \label{1.2}
j = \left[
\begin{array}{cc}
I_{m_1} & 0 \\ 0 & -I_{m_2}
\end{array}
\right], \hspace{1em} V= \left[\begin{array}{cc}
0&v\\v^{*}&0\end{array}\right],
 \end{equation}
where  $I_{m_k}$ is the $m_k \times
m_k$ identity
matrix and $v(x)$ is an $m_1 \times m_2$ matrix function,
which is called the potential of system. 
Dirac-type systems   are very well-known in
mathematics and
applications (see, for instance, books
\cite{AD, CS, K, LS, SaL3},
recent publications \cite{AGKLS1, AGKLS2, AD0,  CG1, CG2, FKS1, FKS2, SaL4},
and numerous references therein).
 The
name ZS-AKNS
is caused by the fact that system (\ref{1.1}) is an
auxiliary
linear system for many important nonlinear integrable
wave
equations and as such it was studied, for instance, in
\cite{AKNS, AS, FT, GSS, SaA7, ZS}.
The Weyl and spectral theory of self-adjoint Dirac systems, where $m_1=m_2$,
was treated, for instance, in \cite{AD0, Kre1, CG1, GKS6, LS, SaA3,  SaL3}
(see also various references therein). The "non-classical" Weyl theory for the
equally important  case $m_1\not= m_2$ 
and related questions are the subject of this paper.  

In Section \ref{DP} we treat the direct problem for the general-type Dirac system,
that is, system \eqref{1.1}, where the potential $v$ is locally summable.
A definition of the non-expansive generalized Weyl function is given,
its existence and uniqueness are proved, and some basic properties are
studied.

In Section \ref{dies} we consider Dirac systems with  so called
{\it generalized pseudo-exponential} potentials (see Definition \ref{Dnpsex}).
Direct and inverse problems for such systems are solved there explicitly.

As usual, $\BN$ stands for the set of natural numbers, $\BR$ stands for the real axis,
$\BC$ stands for the complex plain, and
$\BC_+$ for the open upper
semi-plane. If $a\in \BC$, then $\ov a$ is its complex conjugate.  The notation $\im$
is used for image.
An $m_2 \times m_1$ matrix $\a$ is said to be non-expansive, 
if $\a^*\a \leq I_{m_1}$ (or, equivalently, if $\a\a^* \leq I_{m_2}$).

We put $m_1+m_2=:m$. The fundamental solution of system
\eqref{1.1} is denoted by $u(x,z)$, and this solution is normalized by the condition
\begin{align}&      \label{1.3}
u(0,z)=I_m.
\end{align}
\section{Direct problem} \label{DP}
\setcounter{equation}{0}
We consider Dirac system \eqref{1.1} on the semi-axis $x \in [0,\, \infty)$
and assume that $v$ is measurable and locally summable, that is, summable
on all the finite intervals.
In a way, which is similar, for instance, to the non-classical problem treated in \cite{SaA21}
we shall use M\"obius transformations and matrix balls to solve the direct problem for Dirac
system. 

Introduce a class  of nonsingular $m \times m_1$ matrix functions 
$\clp(z)$ with property-$j$, which are an immediate analog of the classical pairs
of parameter matrix functions. Namely, the matrix functions 
$\clp(z)$ are meromorphic in $\BC_+$ and satisfy
(excluding, possibly, a discrete set of points)
the following relations
\begin{align}\label{2.1}&
\clp(z)^*\clp(z) >0, \quad \clp(z)^* j \clp(z) \geq 0 \quad (z\in \BC_+).
\end{align}
\begin{Dn} \label{set}
The set $\cln(x,z)$ of M\"obius transformations is the set of values at $x, \,z$ 
of matrix functions
\begin{align}\label{2.2}&
\vp(x,z,\clp)=\begin{bmatrix}
0 &I_{m_2}
\end{bmatrix}u(x,z)^{-1}\clp(z)\Big(\begin{bmatrix}
I_{m_1} & 0
\end{bmatrix}u(x,z)^{-1}\clp(z)\Big)^{-1},
\end{align}
where $\clp(z)$ are nonsingular  matrix functions 
 with property-$j$. 
 \end{Dn}
 \begin{Pn} \label{PnW1} Let Dirac system \eqref{1.1} on $[0, \, \infty)$
 be given and assume that $v$ is  locally summable.
 Then the sets $\cln(x,z)$
 are well-defined. There is a unique matrix function
 $\vp(z)$ in $\BC_+$ such that
\begin{align}&      \label{2.3}
\vp(z)=\bigcap_{x<\infty}\cln(x,z).
\end{align} 
This function is analytic and non-expansive.
 \end{Pn}
\begin{proof}.
It is immediate from \eqref{1.1} that
\begin{align}&      \label{2.4}
\frac{d}{dx}\big(u(x,z)^*ju(x,z)\big)=\I (z-\ov z)u(x,z)^*u(x,z)<0, \quad z \in \BC_+.
\end{align}
According to \eqref{1.3} and \eqref{2.4} we have
\begin{align}&      \label{2.5}
\mathfrak{A}(x,z)=\{\mathfrak{A}_{ij}(x,z)\}_{i,j=1}^2:= u(x,z)^*ju(x,z)\leq j,
\quad z \in \BC_+,
\end{align}
where $\mathfrak{A}$ is partitioned into four blocks so that $\mathfrak{A}_{ii}$
is an $m_i \times m_i$ matrix function ($i=1,2$).
Inequality \eqref{2.5} yields
\begin{align}&      \label{2.6}
\big(u(x,z)^*\big)^{-1}ju(x,z)^{-1}\geq j.
\end{align}
Thus, we get 
\begin{align}&      \label{2.7}
\det \Big(\begin{bmatrix}
I_{m_1} & 0
\end{bmatrix}u(x,z)^{-1}\clp(z)\Big)\not= 0,
\end{align}
and so $\cln$ is well-defined via \eqref{2.2}. Indeed, if \eqref{2.7} does not
hold, there is a vector $f \in \BC^{m_1}$ such that
\begin{align}&      \label{2.8}
\begin{bmatrix}
I_{m_1} & 0
\end{bmatrix}j u(x,z)^{-1}\clp(z)f=\begin{bmatrix} I_{m_1} & 0
\end{bmatrix}u(x,z)^{-1}\clp(z)f= 0, \quad f\not=0.
\end{align}
By \eqref{2.1} and \eqref{2.6} the subspace $\im \big(u(x,z)^{-1}\clp(z)\big)$
is a maximal $j$-non-negative subspace. Clearly $\im \big(\begin{bmatrix}
I_{m_1} & 0
\end{bmatrix}^*\big)$ is a maximal $j$-nonnegative subspace too. Therefore
\eqref{2.8} implies $u(x,z)^{-1}\clp(z)f \in \im \big(\begin{bmatrix}
I_{m_1} & 0
\end{bmatrix}^*\big)$. But then it follows from
$$
\begin{bmatrix} I_{m_1} & 0
\end{bmatrix}u(x,z)^{-1}\clp(z)f= 0
$$
that $f=0$, which contradicts the inequality in \eqref{2.8}.

Next, rewrite \eqref{2.2} in the equivalent form
\begin{align}&      \label{2.9}
\begin{bmatrix}
I_{m_1} \\ \vp(x,z,\clp)
\end{bmatrix}=u(x,z)^{-1}\clp(z)\Big(\begin{bmatrix}
I_{m_1} & 0
\end{bmatrix}u(x,z)^{-1}\clp(z)\Big)^{-1}.
\end{align}
In view of \eqref{2.1}, \eqref{2.9}, and of the definition of 
$\mathfrak{A}$ in \eqref{2.5},
formula  
\begin{align}&      \label{2.10}
\wh \vp(z) \in \cln(x,z)
\end{align}
is equivalent to
\begin{align}&      \label{2.11}
\begin{bmatrix}
I_{m_1} & \wh \vp(z)^*
\end{bmatrix}\mathfrak{A}(x,z)\begin{bmatrix}
I_{m_1} \\ \wh \vp(z)
\end{bmatrix} \geq 0.
\end{align}
In a standard way, using formula \eqref{2.4} and the equivalence of \eqref{2.10} and \eqref{2.11}, we get
\begin{align}&      \label{2.16}
\cln(x_1,z) \subset \cln(x_2,z)  \quad {\mathrm{for}} \quad x_1>x_2.
\end{align}
Moreover, \eqref{2.11} at $x=0$ means that
\begin{align}&      \label{2.17}
\cln(0,z) =\{\wh \vp(z): \, \wh \vp(z)^* \wh \vp(z)\leq I_{m_1}\}.
\end{align}
By Montel's theorem, formulas \eqref{2.16} and \eqref{2.17}
imply the existence of an analytic and non-expansive matrix function
$\vp(z)$ such that
\begin{align}&      \label{2.18}
\vp(z)\in \bigcap_{x<\infty}\cln(x,z).
\end{align} 
Indeed, because of \eqref{2.16} and \eqref{2.17} we see that the set 
of functions $\vp(x,z,\clp)$ of the form \eqref{2.2} is uniformly bounded
in $\BC_+$. So, Montel's theorem is applicable and there is an analytic
matrix function, which we denote by $\vp_{\infty}(z)$ and which is a uniform limit
of some sequence
\begin{align}&      \label{M1}
\vp_{\infty}(z)=\lim_{i \to \infty} \vp(x_i,z,\clp_i) \quad (i \in \BN, \quad x_i \uparrow, \quad \lim_{i \to \infty}x_i=\infty)
\end{align} 
on all the bounded and closed subsets of $\BC_+$. Since $x_i \uparrow$
and equalities \eqref{2.9} and \eqref{2.16} hold, it follows that the matrix functions
\begin{align}&      \nn
\clp_{ij}(z):=u(x_i,z)\begin{bmatrix}
I_{m_1} \\ \vp(x_j,z,\clp_j) 
\end{bmatrix} \quad (j \geq i)
\end{align} 
satisfy relations \eqref{2.1}. Therefore, using \eqref{M1} we derive that \eqref{2.1} holds for
\begin{align}&      \nn
\clp_{i, \infty}(z):=u(x_i,z)\begin{bmatrix}
I_{m_1} \\ \vp_{\infty}(z) 
\end{bmatrix},
\end{align} 
which, in its turn, implies that 
\begin{align}&      \label{M2}
\vp_{\infty}(z)\in\cln(x_i,z).
\end{align} 
Since \eqref{M2} holds for all $i\in \BN$,
we see that \eqref{2.18} is true for $\vp(z)=\vp_{\infty}(z)$.

Now, let us show  that $\cln$ is a matrix ball.  It follows from
\eqref{2.5} that $\mf{A}_{22}<-I_{m_2}$. Moreover, it follows 
from \eqref{2.4} and \eqref{2.5} that $\frac{d}{dx}\mf{A} \leq \I(\ov z -z)\mf{A}\leq
 \I(\ov z -z)j$.
Taking into account the inequalities above, we derive
\begin{align}&      \label{2.12}
- \mathfrak{A}_{22}(x,z) \geq \big(1+\I(\ov z -z)x\big)I_{m_2}.
\end{align}
Note also that \eqref{2.5} implies $\mf{A}(x,z)^{-1} \geq j$ 
for $z\in \BC_+$ (see \cite{P}). Thus, we get
\begin{align}&      \label{2.13}
\big(\mathfrak{A}^{-1}\big)_{11}=\big(\mathfrak{A}_{11}-
\mf{A}_{12}\mf{A}_{22}^{-1}\mf{A}_{21}\big)^{-1} \geq I_{m_1}.
\end{align}
As $- \mf{A}_{22}>0$ the square root $\Up=\big(- \mf{A}_{22}\big)^{1/2}$
is well-defined and we rewrite \eqref{2.11} in the form
\begin{align}&      \nn
\mathfrak{A}_{11}-
\mf{A}_{12}\mf{A}_{22}^{-1}\mf{A}_{21}-\big(\wh \vp^*\Up-\mf{A}_{12}\Up^{-1}\big)
\big(\Up \wh \vp-\Up^{-1}\mf{A}_{21}\big) \geq 0, \end{align}
where $\mf{A}_{12}=\mf{A}_{21}^*$. Equivalently, we have
\begin{align}&      \label{2.14}
\wh \vp =\rho_l \om \rho_r-\mf{A}_{22}^{-1}\mf{A}_{21}, \quad \om^*\om \leq I_{m_2},
\\ &      \label{2.15}
 \rho_l:=\Up^{-1}=\big(- \mf{A}_{22}\big)^{-1/2},
\quad \rho_r:=(\mathfrak{A}_{11}-
\mf{A}_{12}\mf{A}_{22}^{-1}\mf{A}_{21}\big)^{1/2}.
\end{align}
Here $\om$ is an $m_2 \times m_1$ matrix function.
Since \eqref{2.10} is equivalent to \eqref{2.14}, the sets
$\cln(x,z)$ (where the values of $x$ and $z$ are fixed) are matrix balls,
indeed.
According to \eqref{2.12}, \eqref{2.13}, and \eqref{2.15}
the next formula holds:
\begin{align}&      \label{2.19}
\rho_l(x,z) \to 0 \quad (x \to \infty), \quad \rho_r(x,z) \leq I_{m_1}.
\end{align}

Finally, relations \eqref{2.18}, \eqref{2.14}, and \eqref{2.19}
imply \eqref{2.3}.

\end{proof}
In view of Proposition \ref{PnW1}
we define the Weyl function similar to the canonical systems 
case \cite{SaL3}. 
\begin{Dn} \label{defWeyl} The Weyl-Titchmarsh (or simply Weyl) function of Dirac system \eqref{1.1} on $[0, \, \infty)$,
where  potential $v$ is locally summable, is the function $\vp$
given by \eqref{2.3}.
\end{Dn}
By Proposition  \ref{PnW1} the Weyl-Titchmarsh function always exists.
Clearly, it is unique.
 
\begin{Cy} \label{CyW2} Let the conditions of Proposition \ref{PnW1} hold.
Then the Weyl function is the unique function, which satisfies the inequality
\begin{align}&      \label{2.20}
\int_0^{\infty}
\begin{bmatrix}
I_{m_1} & \vp(z)^*
\end{bmatrix}
u(x,z)^*u(x,z)
\begin{bmatrix}
I_{m_1} \\ \vp(z)
\end{bmatrix}dx< \infty .
\end{align}
\end{Cy}
\begin{proof}.  According to the equalities in \eqref{2.4} and \eqref{2.5}
and to the inequality \eqref{2.11} we derive
\begin{align}&      \label{2.21}
\int_0^{r}
\begin{bmatrix}
I_{m_1} & \vp(z)^*
\end{bmatrix}
u(x,z)^*u(x,z)
\begin{bmatrix}
I_{m_1} \\ \vp(z)
\end{bmatrix}dx =\frac{\I}{z - \ov z}\begin{bmatrix}
I_{m_1} & \vp(z)^*
\end{bmatrix}
\\ \nn &
\times \big(
\mf{A}(0,z)-\mf{A}(r,z)\big)
\begin{bmatrix}
I_{m_1} \\ \vp(z)
\end{bmatrix}\leq \frac{\I}{z - \ov z}\begin{bmatrix}
I_{m_1} & \vp(z)^*
\end{bmatrix}
\mf{A}(0,z)\begin{bmatrix}
I_{m_1} \\ \vp(z)
\end{bmatrix}.
\end{align}
Inequality \eqref{2.20} is immediate from \eqref{2.21}. Moreover,
as $u^*u\geq -\mf{A}$, the inequality \eqref{2.12} yields
\begin{align}&      \label{2.22}
\int_0^{r}
\begin{bmatrix}0 &
I_{m_2} 
\end{bmatrix}
u(x,z)^*u(x,z)
\begin{bmatrix}
0 \\ I_{m_2} 
\end{bmatrix}dx \geq 
r I_{m_2} .
\end{align}
By \eqref{2.22}, the function satisfying \eqref{2.20} is unique.
\end{proof}
\begin{Rk}\label{RkDefW} In view of Corollary \ref{CyW2}, inequality
\eqref{2.20} can be used as an equivalent definition of the Weyl function.
Definition of the form \eqref{2.20} is a more classical one and deals with
solutions of \eqref{1.1} which belong to $L^2(0, \, \infty)$. Compare with definitions
of Weyl-Titchmarsh or $M$-functions for discrete and continuous
systems in \cite{CG2, LS, Mar, SaA1, SaA2, SaL3, T1, T2} (see also references therein).
\end{Rk}
Our last proposition in this section is dedicated to a property of the Weyl
function,  the analog of which may be used as a definition of  generalized
Weyl functions in more complicated non-self-adjoint cases (see, e.g.,
\cite{FKS1, SaA2, SaA21}). 
\begin{Pn}\label{PnGW} Let Dirac system \eqref{1.1} on $[0, \, \infty)$
 be given, and assume that $v$ is  locally summable.
 Then,  the following inequality
\begin{align}&      \label{6.37}
\sup_{x \leq l, \, z\in \BC_+}\left\|\E^{-\I xz}u(x,z)\begin{bmatrix}
I_{m_1} \\ \vp(z)
\end{bmatrix}\right\|<\infty
\end{align}
holds on  any finite interval $[0, \, l]$
for the Weyl function $\vp$ of this system. 
\end{Pn}
\begin{proof}. We fix some $l$. Now, choose $x$ such that $0<x\leq l<\infty$. Because of \eqref{2.3} the Weyl function $\vp$ admits representations \eqref{2.2} 
(i.e., $\vp(z)=\vp(x,z, \clp)$). Hence, we can use
\eqref{2.1} and \eqref{2.9} to get
\begin{align}&      \label{6.38}
\Psi(x,z)^* j \Psi(x,z) \geq 0, \quad
\Psi(x,z):=\E^{-\I xz}u(x,z)\begin{bmatrix}
I_{m_1} \\ \vp(z)
\end{bmatrix}.
\end{align}
On the other hand, equation \eqref{1.1} and definition of $\Psi$
in \eqref{6.38} imply that
\begin{align}&      \nn
\frac{d}{dx}
\Big(\E^{-2xM}\Psi(x,z)^*(I_m+ j) \Psi(x,z) \Big)
=\E^{-2xM}\Psi(x,z)^*
\\ &      \label{6.39}
\times  \Big(\I \big((I_m+ j)jV-Vj(I_m+ j)\big)-2M(I_m+ j)\Big) \Psi(x,z)
\\ & \nn =2\E^{-2xM}\Psi(x,z)^*\begin{bmatrix}
-2 M I_{m_1} & \I v(x)\\ - \I v(x)^* &0
\end{bmatrix}\Psi(x,z), \quad
M:=\sup_{x<l}\|V(x)\|.
\end{align}
Using \eqref{6.38} and \eqref{6.39} we derive
\begin{align}&      \nn
\frac{d}{dx}
\Big(\E^{-2xM}\Psi(x,z)^*(I_m+ j) \Psi(x,z) \Big) \leq
2\E^{-2xM}\Psi(x,z)^*
\\ & \label{6.40}
\times \left(\begin{bmatrix}
 & \I v(x)\\ - \I v(x)^* &0
\end{bmatrix}-MI_m \right)\Psi(x,z)\leq 0.
\end{align}
Finally, inequalities  \eqref{6.38} and \eqref{6.40}  lead us to
\begin{align}&      \label{6.41}
\Psi(x,z)^* \Psi(x,z) \leq
\Psi(x,z)^*(I_m+ j) \Psi(x,z)  \leq
2\E^{2xM}I_{m_1},
\end{align}
and \eqref{6.37} follows.
\end{proof}
\section{Direct and inverse problem: explicit solutions} \label{dies}
\setcounter{equation}{0}
Various versions of B\"acklund-Darboux transformations are actively
used to construct explicit solutions of linear and integrable nonlinear
equations (see, e.g., \cite{Ci, GeH, Gu, Mar2, MS, ALS10, ZM} and numerous
references therein). For the spectral and scattering results that follow
from B\"acklund-Darboux transformations and related commutation and factorization methods see, for instance, publications \cite{cr, de, FKS0, fg, gt, GKS6, KoSaTe,
kr57, SaA7}. Here we will give explicit solutions of our  direct and inverse problems using the GBDT version of the B\"acklund-Darboux transformation
(see \cite{FKRS, FKS0, GKS6, ALS94, SaA7, ALS10} and references therein).

To obtain explicit solutions, we consider
$m_1 \times m_2$
potentials $v$ of the form
\begin{align} \label{7.1}&
v(x)=-2i \vt _{1}^{\, *}\E^{ix \alpha^{*}} \S
(x)^{-1}\E^{ix
\alpha }\vt_2,
\end{align}
where some  $n\in \BN$ is fixed and
the $n \times n$ matrix function $\S$ 
is given by
the formula
\begin{align}   \label{7.2}&
\S(x)=\S_0+ \int_{0}^{x} \Lambda(t) \Lambda (t)^{*}dt
\quad(\S_0>0),
\quad
\Lambda (x)= \begin{bmatrix}  \E^{-ix \alpha } \vt_{1} 
& \E^{ix
\alpha } \vt_2 \end{bmatrix}.
\end{align}
Here $\a$, $\vt_1$, and $\vt_2$ are $n \times n$, $n \times
m_1$, and $n
\times m_2$ parameter matrices,  and the following
matrix identity holds:
\begin{align} \label{7.3}&
\a \S_0-\S_0 \a^*=\I (\vt_1 \vt_1^*- \vt_2 \vt_2^*).
\end{align}
Clearly, $\S(x)$ is invertible for $x\geq 0$ and the potential $v$ in \eqref{7.1}
is well-defined. 
\begin{Dn}\label{Dnpsex}
The $m_1 \times m_2$
potentials $v$ of the form \eqref{7.1}, where relations \eqref{7.2} and \eqref{7.3}
hold, are called generalized pseudo-exponential potentials. It is said that
$v$ is generated by the parameter matrices $\a$, $\S_0$, $\vt_1$, and $\vt_2$.
\end{Dn}
According to \cite[Theorem 3]{ALS94}
(see also \cite{FKS0}), the fundamental solution $u$ of system \eqref{1.1}, where
$V$ is given by \eqref{1.2}, $v$   is a generalized pseudo-exponential potential, and $u$ is
normalized by \eqref{1.3}, admits representation
\begin{align} \label{7.4}&
u(x, z)=w_{\a}(x, z)\E^{ix z j}w_{\a}(0, z)^{-1}.
\end{align}
Here we have
\begin{align} \label{7.5}&
w_{\a}(x, z):=I_m+ij \Lam(x)^*\S(x)^{-1}(z I_n- \a)^{-1}
\Lam(x).
\end{align}

Note that the case  $m_1=m_2$ (i.e., the case of the 
pseudo-exponential potentials) was treated in greater detail in \cite{GKS6}
(see \cite{GKS6} and references therein for the term {\it pseudo-exponential},
itself, too).  Formulas (\ref{7.2}) and (\ref{7.3}) yield
\begin{equation} \label{7.6}
\a \S(x)-\S(x) \a^*=i \Lam (x)j \Lam (x)^*.
\end{equation}
Identity \eqref{7.6}, in its turn, implies that $w_{\a}(z)$  
is a transfer matrix function
in Lev Sakhnovich form \cite{SaL1, SaL2, SaL2', SaL3}.
However, $w_{\a}(x, z)$  possesses an additional variable
$x$ and the way, in which this matrix function depends on $x$,
is essential. 

From \cite[formula (2.9)]{FKS0}, where $W_{11}$ and $W_{21}$
are $m_1 \times m_1$ and $m_2 \times m_1$ blocks of $w_{\a}(0, z)$:
\begin{align} \label{7.7}&
w_{\a}(0, z)=:\{W_{ij}(z)\}_{i,j=1}^2,
\end{align}
we see that
\begin{align} \label{7.8}&
W_{21}(z)W_{11}(z)^{-1}= - i\vt_2^*\S_0^{-1}(z I_n
- \t)^{-1}
\vt_1 , \quad \t:=\a-i\vt_1\vt_1^*\S_0^{-1}.
\end{align}
We note that  \cite[formulas (2.6)  and (2.7)]{FKS0} imply that $W_{11}(z)$ 
is always well-defined and invertible for $z\not\in \s(\a)\cup\s(\t)$, where
$\s$ denotes spectrum.

Relations \eqref{7.4}, \eqref{7.7}, and \eqref{7.8} are basic to solve the direct
problem for Dirac systems with the generalized pseudo-exponential potentials
\eqref{7.1}.
\begin{Tm}\label{TmDPes} Let Dirac system \eqref{1.1} on $[0, \, \infty)$
 be given and assume that $v$ is  a generalized pseudo-exponential potential,
 which is generated by the matrices $\a$, $\S_0$, $\vt_1$, and $\vt_2$.
 Then the Weyl function $\vp$ of  system \eqref{1.1} has the form:
\begin{align} \label{7.9}&
\vp(z)=- i\vt_2^*\S_0^{-1}(z I_n
- \t)^{-1}
\vt_1 , \quad \t=\a-i\vt_1\vt_1^*\S_0^{-1}.
\end{align}
\end{Tm}
\begin{proof}. We compare \eqref{7.8} and \eqref{7.9} to see that
\begin{align} \label{7.10}&
\vp(z)=W_{21}(z)W_{11}(z)^{-1}.
\end{align}
Because of \eqref{7.4}, \eqref{7.7}, and \eqref{7.10} we have
\begin{align} \label{7.11}&
u(x,z)
\begin{bmatrix}
I_{m_1} \\ \vp(z)
\end{bmatrix}=\E^{ixz}w_{\a}(x, z)\begin{bmatrix}
I_{m_1} \\ 0
\end{bmatrix}W_{11}(z)^{-1}.
\end{align}
To consider the matrix function $\Lam^*\S^{-1}$, which appears in the definition \eqref{7.5} of $w_{\a}$,
we derive from \eqref{7.2} that
\begin{align} \label{7.12}&
\S(x)^{-1}\Lam(x)\Lam(x)^*\S(x)^{-1}= -\frac{d}{dx}\S(x)^{-1}.
\end{align}
It is immediate also from \eqref{7.2} that $\S(x)>0$. 
Therefore, using  \eqref{7.12} we get 
\begin{align} \label{7.13}&
\int_0^{\infty} \S(t)^{-1}\Lam(t)\Lam(t)^*\S(t)^{-1}dt\leq \S_0^{-1}.
\end{align}
Furthermore, the last equality in \eqref{7.2} implies that
\begin{align} \label{7.14}&
\sup_{\Im z>\|\a\|+\ve }\|\E^{ixz}\Lam(x)\|<M_{\ve} \quad (\ve>0).
\end{align}
It follows from \eqref{7.5}, \eqref{7.13}, and \eqref{7.14} that
the entries of the right-hand side of \eqref{7.11} are well-defined and
uniformly bounded in the $L^2(0,\, \infty)$ norm with respect to $x$
for all $z$ such that $\Im z \geq \max \big(\|\a\|, \,  \|\t\|\big)+\ve$ and $\ve>0$.
Hence, taking into account \eqref{7.11} we see that \eqref{2.20} holds for $z$
from the  mentioned above domain. So, according to the
uniqueness statement in
Corollary \ref{CyW2},
$\vp(z)$ of the form \eqref{7.9} coincides with the Weyl function in that domain. Since the Weyl function
is analytic in $\BC_+$, the matrix function $\vp$ coincides with it in $\BC_+$
(i.e., $\vp$ is the Weyl function, indeed).
\end{proof}
For the case that $v$ is a  
generalized pseudo-exponential potential, where $\S_0>0$, our Weyl function 
coincides with the reflection coefficient from \cite{FKS0} (see 
\cite[Theorem 3.3]{FKS0}). Hence, the solution of our inverse problem
can be considered as a particular case of  the solution of the inverse problem
from \cite[Theorem 4.1]{FKS0}, where the singular case $\S_0\not>0$ was
treated too.

Before we formulate the procedure to solve inverse problem, some 
results  on rational matrix functions and notions from system and  control theories are required (see, e.g., \cite{KFA, LR}). Let $\vp(z)$ be a strictly proper rational matrix function, that is, such a rational matrix function that
\begin{align} \label{7.15}&
\lim_{z\to \infty} \vp(z)=0.
\end{align}
Then $\vp$ admits representations (also called realizations):
\begin{align} \label{7.16}&
 \vp(z)=C_N(z I_N-\cla_N)^{-1}B_N,
\end{align}
where $C_N$, $\cla_N$, and $B_N$ are $m_2 \times N$,
$N \times N$, and $N \times m_1 $, respectively, matrices. Here $N\in \BN$,
and $m_1$ ($m_2$) denotes the number of columns (rows) of $\vp$.
\begin{Dn}\label{DnST}
The minimal possible value of $N$ in realizations \eqref{7.16} is called
the McMillan degree of $\vp$, and we denote this value by $n$.
Realizations \eqref{7.16}, where $N=n$, are called the minimal realizations.
\end{Dn}
From \cite[Theorems 21.1.3, 21.2.1]{LR} we easily see that for  a minimal
realization
\begin{align} \label{7.17}&
 \vp(z)=C(z I_n-\cla)^{-1}B
\end{align}
of a matrix $\vp$, which is non-expansive on $\BR$ and has no poles
in $\BC_+$, there is a positive solution $X>0$ of the Riccati equation 
\begin{align} \label{7.18}&
XC^*CX+\I(X\cla^*-\cla X)+BB^*=0.
\end{align}
Furthermore, all the hermitian
solutions of \eqref{7.18} are positive.
\begin{Tm}\label{TmIpes} Let $\vp(z)$ 
be a strictly proper rational matrix function, which  is non-expansive on $\BR$ and has no poles in $\BC_+$.  Assume that \eqref{7.17} is its minimal realization
and that $X>0$ is a solution of \eqref{7.18}. 

Then $\vp(z)$ is the Weyl function of Dirac system, the potential of which
is given by \eqref{7.1} and \eqref{7.2}, where
\begin{align} \label{7.19}&
\a=\cla +\I BB^*X^{-1}, \quad \S_0=X, \quad \vt_1=B, \quad \vt_2=-\I X C^*.
\end{align}
This solution of the inverse problem is unique in the class of Dirac systems
with locally bounded potentials.
\end{Tm}
\begin{proof}.  From \eqref{7.19} we see that
\[
\a\S_0-\S_0\a^*=\cla X-X\cla^*+2iBB^*, \quad
\I(\vt_1\vt_1^*-\vt_2\vt_2^*)=\I BB^*-\I XC^*CX,
\]
and so \eqref{7.3} is equivalent to \eqref{7.18}. Since \eqref{7.3} holds,
we apply Theorem \ref{TmDPes}.   Theorem \ref{TmDPes} implies that the Weyl function of the Dirac system, where $v$ is given by \eqref{7.1},
has the form \eqref{7.9}.
Next, we substitute \eqref{7.19} into \eqref{7.9}, to derive that the right-hand
sides of \eqref{7.17} and the first equality in \eqref{7.9} coincide.
In other words, the Weyl function of our system admits representation
\eqref{7.17}.

Finally, the uniqueness of the solution of the inverse problem
follows from Theorem \cite[Theorem 4.1]{FKRS3}.
\end{proof}
We note that the corresponding uniqueness result in \cite{FKS0}
was proved only for the class of systems with the generalized
pseudo-exponential potentials.

Because of the second equality in  \eqref{7.8} and identity
\eqref{7.3}, the matrix $\t$ satisfies another identity:
$\t \S_0-\S_0 \t^*=-\I (\vt_1 \vt_1^*+ \vt_2 \vt_2^*)$, that is,
\begin{align} \label{7.22}&
 \S_0^{-1}\t -\t^* \S_0^{-1}=-\I \S_0^{-1} (\vt_1 \vt_1^*+ \vt_2 \vt_2^*)\S_0^{-1}.
\end{align}
If $f\not=0$ is an eigenvector of $\t$ (i.e., $\t f=\la f$), identity \eqref{7.22}
implies that
\begin{align}  \label{7.23}&
(\la - \ov{\la})f^*\S_0^{-1}f=-\I f^*\S_0^{-1} (\vt_1 \vt_1^*+ \vt_2 \vt_2^*)\S_0^{-1}f.
\end{align}
Since $\S_0>0$, we derive from \eqref{7.23} that
\begin{align}  \label{7.24}&
\s(\t) \subset \BC_- \cup \BR .
\end{align}

Real eigenvalues of $\t$ play a special role in the spectral
theory of  an operator, which corresponds to the
Dirac system with a generalized pseudo-exponential potential
(see, e.g., \cite{GKS5} for the case of square potentials).
In our case the operator $\clh$  corresponding to the
Dirac system is defined in a way, which is similar
to  the definition from \cite{GKS5}, but the initial condition is
quite different. Namely, we determine $\clh$ 
by the differential expression
\begin{align} \label{7.20}&
\clh_{de} y=-\I j \frac{d}{dx}y-Vy,
\end{align}
and by its domain $\cld(\clh)$, which consists of all locally absolutely
continuous $\BC^m$-valued functions $y$ in $L^2_m(0, \, \infty)$, 
such that 
\begin{align} \label{7.21}&
\clh_{de} y\in L^2_m(0, \, \infty), \qquad y(0)=0.
\end{align}
\begin{Pn}\label{PnEV} Let the conditions of Theorem \ref{TmDPes}
hold,  let $\t$ be given by the second relation in \eqref{7.8}, and let
$\la$ be an eigenvalue of $\t$:
\begin{align} \label{7.25}&
\t f=\la f, \quad f \not= 0, \quad \la \in \BR .
\end{align}
Then, the matrix function
\begin{align} \label{7.26}&
g(x):=j\Lam(x)^*\S(x)^{-1}f
\end{align}
is a bounded state of $\clh$ and $\clh g=\la g$.
\end{Pn}
\begin{proof}.  First, we show that formulas \eqref{7.23} and \eqref{7.25} yield
\begin{align} \label{7.27}&
 \vt_1^*\S_0^{-1}f=0, \quad    \vt_2^*\S_0^{-1}f=0, \quad          \a f=\la f.
 \end{align}
Indeed, the first two equalities in \eqref{7.27} easily follow from
\eqref{7.23} for the case that $\la =\ov \la$. The equality $ \a f=\la f$
is immediate from $\t f = \la f$, definition of $\t$ in \eqref{7.8}, and equality
$\vt_1^*\S_0^{-1}f=0$.

Next, we show that
\begin{align} \label{7.28}&
\big(j\Lam^*\S^{-1}\big)^{\prime}=\I j^2 \Lam^*\S^{-1}\a+
\big(\Lam^*\S^{-1}\Lam - j \Lam^*\S^{-1}\Lam j\big)j \Lam^*\S^{-1}.
\end{align}
Formula \eqref{7.28} follows from  a general GBDT formula
\cite[(3.14)]{ALS10} and also from its Dirac system subcase 
\cite[(2.13)]{ALS10}, but it will be convenient to prove \eqref{7.28} directly.
We note that  formula \eqref{7.2} implies
\begin{align} \label{7.29}&
\Lam^{\prime}= -\I \a \Lam j, \qquad \S^{\prime}=\Lam\Lam^*,
\end{align}
and formula \eqref{7.6} can be rewritten as
\begin{align} \label{7.30}&
\a^* \S^{-1}=\S^{-1}\a- \I \S^{-1}\Lam j \Lam^*\S^{-1}.
\end{align}
Since $j^2=I_m$, using \eqref{7.29} and \eqref{7.30} we obtain \eqref{7.28}.

Now, partitioning $\Lam$ into two blocks and using \eqref{7.1} and \eqref{7.2}, we see that
\begin{align} \label{7.31}&
v(x)=-2\I \Lam_1(x)^*\S(x)^{-1}\Lam_2(x), \quad \Lam=: \begin{bmatrix} \Lam_1 & \Lam_2 \end{bmatrix}.
\end{align}
In view of \eqref{1.2} and  \eqref{7.31} we have
\begin{align} \label{7.32}&
\Lam^*\S^{-1}\Lam - j \Lam^*\S^{-1}\Lam j=\I j V.
\end{align}
Applying both sides of \eqref{7.28} to $f$ and taking into account the last equality in \eqref{7.27} and relation \eqref{7.32},
we derive
\begin{align} \label{7.33}&
\big(j\Lam(x)^*\S(x)^{-1}f\big)^{\prime}=\I \la j^2 \Lam(x)^*\S(x)^{-1}f+
\I jV(x)j \Lam(x)^*\S(x)^{-1}f.
\end{align}
Because of  \eqref{7.20} and  \eqref{7.26}, we can rewrite \eqref{7.33} as
\begin{align} \label{7.34}&
\clh_{de}g=\la g,
\end{align}
and it remains to show that $g\in \cld(\clh)$, that is, that $g \in L^2_m(0,\, \infty)$ and \eqref{7.21}
holds for $y=g$. From  \eqref{7.13} and  \eqref{7.34} we see that $g, \, \clh_{de}g \in L^2_m(0,\, \infty)$.
Finally, the initial condition 
\begin{align} \label{7.35}&
g(0)=\begin{bmatrix} \vt_1 & - \vt_2 \end{bmatrix}^*\S_0^{-1}f=0
\end{align}
 is immediate from \eqref{7.2}, \eqref{7.26}, and \eqref{7.27}.
\end{proof}


{\bf Acknowledgement.}
The work of I.Ya. Roitberg was supported by the 
German Research Foundation (DFG) under grant no. KI 760/3-1 and
the work of A.L. Sakhnovich was supported by the Austrian Science Fund (FWF) under Grant  no. Y330.



\end{document}